\documentclass[11pt,english]{article}
\makeatletter
\usepackage{geometry}
\usepackage{float}
\usepackage{graphicx}
\usepackage{amsmath,amsthm}
\usepackage[american]{babel}
\newcommand{\doublespacing}{\let\CS=\@currsize\renewcommand{\baselinestretch}{1.2}\tiny\CS}

\usepackage{hyperref}
\hypersetup{colorlinks=true,linkcolor=blue,citecolor=blue}
\date{}
\begin{document}
\newpage
\doublespacing
\doublespacing
\title {\Large {\bf\it  Optimum Production  through variational principle with the time quadratic demand,  fuzzy time period  and fuzzy integrand}}
\author{{\bf J. N. Roul $^{a},$ K. Maity  $^{b}, $   S. Kar $^{c},$ and M. Maiti $^{d}$}\\
  $^{a}$Patha Bhavana, Visva-Bharati, Santiniketan-731235, West Bengal
\\  India, Email$-$jotin2008@rediffmail.com  \\
  $^{b}$ Department of Mathematics, M G
  Mahavidyalaya, Bhupatinagar-721425\\   
 W.B., India, Email$-$kalipada\_maity@yahoo.co.in\\
 $^{c}$ Department of Mathematics, National Institute of
Technology, Durgapur-713209\\ W.B., India, Email$-$kar\_s\_k@yahoo.com\\
 $^{d}$Department of Applied Mathematics with Oceanology and Computer
 Programming\\
 Vidyasagar University, Midnapore-711201,W.B., India,
Email$-$mmaiti2005@yahoo.co.in}

\maketitle
\begin{abstract}
Here a real life  optimal control problem under fuzzy time period using variational principle is formulated and Solved. The unit production cost is a function of production rate and also dependent on raw material cost, development cost due to durability and wear-tear cost. The holding cost is assumed to be non-linear, dependent on time.  The profit function which consists of  revenue, production cost and holding cost is formulated as a Fuzzy-Final Time and Fixed State System optimal control problem with fuzzy time period. Here production rate is unknown and considered as a control variable and  stock level is taken as a state variable. It  is formulated to optimize the production rate so that total profit is  maximum. The non-linear optimization technique-Generalised Reduced Gradient Method (LINGO 11.0) is used. The optimum results are illustrated both numerically and graphically. 
\end{abstract}
\indent\indent\indent \section*{Submitting our work as a preprint.}
\section{Introduction}
\noindent  Taking on variation of the defined curve and using Euler's necessary conditions and Lagrange's undetermined multipliers, a  mathematical discipline has been applied here. This is known as "Calculus of Variations".  Optimal control theory is an extension of the calculus of variations. \\
Since the development of EOQ model by Harris \cite{harris1913}, a lot of research works has been reported in the literature (cf. Naddor\cite{naddor1966}, Hadly and Whitin \cite{hadleywhitin1975} and others).  In classical inventory models, normally static lot size models are formulated. Due to dynamic manufacturing environment, the static models are not adequate in analyzing the behavior of such systems and in designing the optimal policies for their control. Because of the fact, dynamic models of production inventory systems have been considered by some researchers (cf. Kirk\cite{kirk1970}, Hu and Loulon \cite{huloulon1995},  Hu and Dong \cite{hudong1995}, Khouja \cite{khouja1995}, Giri et al. \cite{giri1996}  and others). In these models, demand and/or production are assumed to be continuous functions of time. During the last two decades, many researchers (cf.  Salameh et al. \cite{salameh2000},  Zequeria et al. \cite{zequeria2004}, Maity and Maiti \cite{maity2007}, Panda et al. \cite{panda2008}, Benjaffar et al. \cite{benjaffar2010}, Sarkar et al.\cite{sarkar2011},  Roul et al. \cite{roul2015} and others) have given considerable attention to the area of inventory of deteriorating/defective/ perishable items, since the life time of an item is not infinite while it is in storage and/or all units can't be produced exactly as per the determined features.  \\ 
Recently Farhadinia \cite{farhadinia2011} developed necessary optimality conditions for fuzzy variational problems, Najariyan and Farahi  \cite{najariyanfarahi2013} reported on optimal control of fuzzy linear
controlled system with fuzzy initial conditions and  Roul et al. \cite{roul2017} a formulated and solved application of  fuzzy variational principle in inventory problems. But, the above mentioned models are developed for the fixed (deterministic) time horizon. Till now, it has been considered that each is exhausted after a fixed time period and this time period is same for all items. This is a hard assumption on the real world situations. Actually, every year a seasonable product does not end at a particular time. In a year, several season products do have also different time periods. There is an inherent uncertainty in these time horizons. This uncertainty can be represented by fuzzy number. Thus the time periods of the items  made out of seasonable products are fuzzy. This is also true in the case of fast moving items also.  \\
The general model is formulated as an optimal control problem over a fuzzy time horizon with profit maximization. This is very simple application on fuzzy time period through optimal control theory. By the theory of Fuzzy  number, fuzzy time period is converted into crisp one by using $\alpha$ measure techniques. Then, the equivalent crisp problem is solved by  Variational Principle, fixed-final time and free-final state system (cf. Naidu \cite{naidu2000})and Generalised Reduced Gradient (GRG) Technique(cf. Gabriel and Ragsdell \cite{gabrielragsdell1977}). The behaviour of optimal production rate, stock and demand with respect to $'\alpha'$ over the time are depicted. The total Maximum profit with respect to different values of $\alpha$  are  pointed. For $\alpha=1.00$, the Max Profit in imprecise environment is equivalent to the result in crisp environment which is established here.
\subsection{Fuzzy algebra:}
\noindent{\bf Fuzzy Number:} A fuzzy number is simply an ordinary number whose precise value is somewhat uncertain. If a fuzzy set is convex and normalized, and its membership function is defined in $R$ and piecewise continuous, it is called as fuzzy number.\\

  Among the various shapes of fuzzy number, triangular fuzzy number is the most popular one. \\
\noindent {\bf Triangular Fuzzy Number (TFN):}\index{Triangular Fuzzy Number} A TFN $\widetilde{A}$ is specified by the triplet
$(a_{1},a_{2},a_{3})$ and is defined by its continuous membership function
 $\mu_{\widetilde{A}}(x): X \,\rightarrow [0,1]$ as follows (cf. Figure 01).
\begin{eqnarray}
\mu_{\widetilde{A}}(x) = \left\{ \begin{array}{ll}
 \displaystyle\frac{x -
a_{1}}{a_{2}-a_{1}} &  \mbox{if} \,\,\,\,a_{1} \leq x \leq a_{2}\\
 \displaystyle \frac{ a_{3}-x}{a_{3}-a_{2}} &
\mbox{if} \,\,\,\,a_{2} \leq x \leq a_{3}\\
\displaystyle \,\,\,\,\,\,\,\, 0  &   \mbox{otherwise} \\
\end{array} \right.\hspace{2 cm}\nonumber
\end{eqnarray}
\begin{figure}
\centering
\includegraphics[width=7.5cm, height=5.0cm]{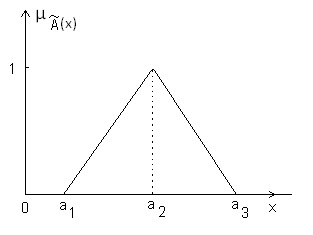}
\caption{Membership function of Triangular Fuzzy Number}
\end{figure}
\noindent {\bf $\alpha$-cut of Fuzzy Number:} If a fuzzy number represented with three points as follows :
$\tilde A = (a_1, a_2, a_3).$   The $\alpha$-cut of fuzzy number $\tilde A $ determines  it crisp one. So $A_{\alpha} = [A^l_{\alpha},A^r _{\alpha}]=[a_1+(a_2-a_1)\alpha, a_3-(a_3-a_2)\alpha] $ for $\alpha \epsilon[0, 1]. $

In practice, we always consider the membership function of the fuzzy number as having bell shape. Thus, we suggest the fuzzy real number $\tilde t$ having membership function (with spread $\sigma$)
\begin{eqnarray}
 \mu(\xi_t) =\left\{\begin{array}{l}
1-\frac{|{\xi_t-t}|}{\sigma}  \mbox{ \hspace {1. cm} if } -\sigma +t \leq  \xi_t \leq \sigma +t\\
 0 \mbox{ \hspace {3.5 cm} otherwise }\end{array}\right.
\nonumber \end{eqnarray}
For final time T which is fuzzy in nature i.e $\tilde T$. It is defined by the following membership function. 
\begin{eqnarray}
 \mu(\xi_T) =\left\{\begin{array}{l}
1-\frac{|{\xi_T-T}|}{\sigma}  \mbox{ \hspace {1. cm} if } -\sigma +T \leq  \xi_T \leq \sigma +T\\
 0 \mbox{ \hspace {3.5 cm} otherwise }\end{array}\right.
 \nonumber\end{eqnarray}
Let  $\alpha=1-\frac{|{\xi_T-T}|}{\sigma}$.\\
Then $-\sigma(1-\alpha)+T \leq  \xi_T  \leq \sigma(1-\alpha)+T $. Therefore $T^{l}(\alpha)=-\sigma(1-\alpha)+T $ and $T^{r}(\alpha)=\sigma(1-\alpha)+T $
\section{Proposed Optimal Control Models }
For  the optimal control problem under finite time horizon, following assumptions and notations are used.
\subsection{Assumptions}
$(i)$ It is a single period production inventory model with
finite time horizon $t\epsilon[0,T]$;\\
$(ii)$ Break-ability rate is unknown;\\
$(iii)$ Shortages are not allowed;\\
$(iv)$ There is no repair or replacement of defective units over whole time period;\\
$(v)$  The development cost of the system increases the durability of the system;\\
$(vi)$ Unit production cost depends on produced-quantity, raw material, wear-tear and development costs;\\
$(vii)$ Deterministic  variational principle is  used for the model.
\subsection{Notations}
{$T :$} total time for the system;\\
{$u(t):$} production rate at time $t$ which is taken as control variable;\\
{$x(t):$} stock level at  time t which is a state variable;\\
{$h(t) :$}=$a+b t( a\,\,\,\, \mbox{and} \,\,\,\,b \,\,\,\,\mbox{are constants })$ holding cost per unit;\\ 
{$L :$} fixed cost like labour, energy, etc;.\\  
{$N :$} cost of technology, design, complexity, resources, etc;.\\  
{$c_{10}:$} constant material cost;\\
{$d_{1}:$} constant demand at initial stage;\\
{$d(t)=d_{1}+d_2 t+d_3 t^2:$} demand function with $d_2$ and $d_3$ are constants.\\
{$c_d(l)=N+L: $} development cost to improve the quality of the product;\\
{$\beta_{10}:$} wear-tear cost for the system.\\
{$p:$} is the selling price of market demand. 

\subsection{Formulation of optimal control  models for dynamic demand with holding cost $h(t)=a+bt$}
\indent In this model,  the differential equation for stock level
{$x(t)$} regarding above system
during a fuzzy time-horizon, $\tilde T $ is
\begin{eqnarray} \frac{d}{dt}\bigg(x(t)\bigg)= u (t)-d(t)
\end{eqnarray}
where $ d(t)=d_1 +d_2 t+d_3 t^2$ \\
 \noindent The unit
production cost is considered as a function of produced-quantity,
raw material cost, wear-tear and development costs(cf. Khouja \cite{khouja1995}). So the total
production cost  is
\begin{eqnarray}
 c_{u} (t)u(t)=\bigg(c_{10}u(t) +{c_{d} (l
)} +\beta _{10} u^2 (t)\bigg)\end{eqnarray} \noindent
\noindent  The total holding cost over the finite time interval
$[0,\,\tilde T]$
 for the stock $x(t)$ is
\begin{eqnarray}
 \int _{0}^{\tilde T} h(t) x (t)\,dt  
 \end{eqnarray}
 where $ h(t)=a+b t $ and $ a,b $ are being constants. \\
 Set-up cost is normally constant with respect to time. But, if dynamic
production rate is considered, some machineries, etc., are to be
set-up and maintained in such a way that the production
system can stand with the pressure of increasing demand.
Thus  a part of set up cost is linearly proportional to production rate and
hence the form is assumed as  $\frac{s_1 +s_2 u(t)}{T}$. Here $s_1$ and $ s_2$  are constants.
Therefore the set up cost $s(t)$ for the model is taken as
\begin{eqnarray}
s(t)=\frac{s_1 +s_2 u(t)}{T}
\end{eqnarray}
 Let $p$  be the selling price of an unit.  Then the revenue from market demand is 
\begin{eqnarray}
\int_{0}^{\tilde T} p d(t) dt
\end{eqnarray}
 Thus the problem reduces to maximization of the profit function $J$
subject to the  constraint satisfying the dynamic
production-demand relation.
\begin{eqnarray}
&&\mbox{Max} J=\int _{0}^{\tilde {T}}\bigg(p d(t) - h(t) x(t)-c_{u }
 (t) u (t)-s(t)\bigg) dt\nonumber\\
&&\mbox{sub\, to}\,\,
 \frac{d}{dt}\bigg(x(t)\bigg)= u (t)-d(t)\\
  &&\mbox{ where \,} u(t)\ge 0\,\,  ,x(t)\ge0\\
\nonumber\end{eqnarray}
\noindent Then the expression (6) is written as
\begin{eqnarray}
 J_{max}\bigg({x}(t),\dot{x}(t),t\bigg)
 &=&\int _{0}^{\tilde T}\bigg(pd(t)-h(t)\, x(t)-\bigg(c_{10}+\frac{ s_2}{T}\bigg)
 \bigg({\dot{x}(t)+d(t)}\bigg) \nonumber\\&-&{ \beta_{10}}\bigg({\dot{x}(t)+d(t)}\bigg)^2 -\bigg({c_{d}(l)}+\frac{ s_1}{T}\bigg)\bigg) dt   \,\, \mbox{ where \,} \,  x(t)\ge0
 \end{eqnarray}
 \noindent The above problem (8) is defined as an optimal
control problem with state variable $x(t)$. Here, (8) contains $u(t)$ implicitly.\\ 
\noindent Using Euler-Lagrange's equation, Fuzzy-Final Time and Fixed State System (i.e, here final
time $\tilde T$
is not specified and $x(t)\equiv x, \dot x(t)\equiv \dot x $, $x(0)=0$ , $ x(T^{l}(\alpha))=0$ for left cut of $\tilde T$; $x(0)=0$, $ x(T^{r}(\alpha))=0$ for right cut of $\tilde T$); we
have

\begin{eqnarray}\frac{\partial
J}{\partial {x}}&-&\frac{d\ }{dt}\bigg(\frac{\partial
J}{\partial \dot{x}}\bigg)=0\\
\nonumber\end{eqnarray}  

\noindent Using equation (9) in (8), we get \\
\begin{eqnarray}
a+bt- 2 \beta_{10} \frac{d\ }{dt} \bigg( \dot{x}(t)+d(t)\bigg)=0
\end{eqnarray} 
Representing $x(t)\equiv x$,
 \begin{eqnarray}
\frac{d^2 x}{dt^2}  =f(t)  \,\,\,\,\mbox{where}\,\,\, f(t)=-2 d_3 t+\frac{ b t}{2\beta_{10}}  +\frac{ a}{2\beta_{10}}-d_2
\end{eqnarray}
This is the second order first degree differential equation of the independent variable t.\\

\noindent Therefore, solving the above we can get  
\begin{eqnarray}
x^l(t,\alpha)&=&A+Bt -(2d_3 -\frac{ b }{2\beta_{10}})\frac{t^3}{6}+(\frac{ a}{2\beta_{10}}-d_2) \frac{t^2}{2}\\
\mbox{and}&&\nonumber\\
x^r(t,\alpha)&=&C+Dt -(2d_3 -\frac{ b }{2\beta_{10}})\frac{t^3}{6}+(\frac{ a}{2\beta_{10}}-d_2) \frac{t^2}{2}
\end{eqnarray}
 Here the initial condition  $x(0)=0$  and final condition
  $x(T^l(\alpha))=0$. Using these conditions, we can get the values of  $A=0$  and $B=(2d_3 -\frac{ b }{2\beta_{10}})\frac{(T^l(\alpha))^2}{6}+(\frac{ a}{2\beta_{10}}-d_2) \frac{T^l(\alpha)}{2}
$. \\
 Similarly,  the initial condition  $x(0)=0$  and final condition
  $x(T^l(\alpha))=0$. Using these conditions, we can get the values of  $C=0$  and $D=(2d_3 -\frac{ b }{2\beta_{10}})\frac{(T^r(\alpha))^2}{6}+(\frac{ a}{2\beta_{10}}-d_2) \frac{T^r(\alpha)}{2}
$. \\
From(1), we have the production rate 
 \begin{eqnarray}
 u^l(t,\alpha)&=&B+ d_1 +\frac{ a}{2\beta_{10}} t+\frac{ b }{2\beta_{10}}\frac{t^2}{2}\\
   \mbox{and}\nonumber\\
 u^r(t, \alpha)&=&D+ d_1 +\frac{ a}{2\beta_{10}} t+\frac{ b }{2\beta_{10}}\frac{t^2}{2}
   \end{eqnarray}
 Then the  profit functional $J$ takes the form $J^l_{max} (\alpha)$
 \begin{eqnarray}
 &=&p(d_1 T^l(\alpha) +d_2 \frac{(T^l(\alpha))^2}{2}+d_3 \frac{(T^l(\alpha))^3}{3})-\bigg(c_d(l)+\frac{s_1}{T^l(\alpha)}\bigg)T^l(\alpha)\nonumber\\
 &-&a\bigg( A T^l(\alpha)+B \frac{(T^l(\alpha))^2}{2} -(2d_3 -\frac{ b }{2\beta_{10}})\frac{(T^l(\alpha))^4}{24}+(\frac{ a}{2\beta_{10}}-d_2) \frac{(T^l(\alpha))^3}{6}
\bigg)\nonumber\\
 &-&b\bigg( A \frac{T^l(\alpha)}{2}+B \frac{(T^l(\alpha))^3}{3} -(2d_3 -\frac{ b }{2\beta_{10}})\frac{(T^l(\alpha))^5}{30}+(\frac{ a}{2\beta_{10}}-d_2) \frac{(T^l(\alpha))^4}{8}\bigg)\nonumber\\
 &-&\bigg(c_{10}+\frac{s_2}{T}\bigg)\bigg( B T^l(\alpha)+ d_1 T^l(\alpha) +\frac{ b }{2\beta_{10}}\frac{(T^l(\alpha))^3}{6}+\frac{ a}{2\beta_{10}} \frac{(T^l(\alpha))^2}{2}\bigg) \nonumber\\
 &-&\beta_{10}\bigg( (B+ d_1)^2 +(\frac{ b }{2\beta_{10}})^2\frac{(T^l(\alpha))^5}{20}+(\frac{ a}{2\beta_{10}})^2 \frac{(T^l(\alpha))^3}{3}+\frac{ b(B+d_1) }{\beta_{10}}\frac{(T^l(\alpha))^3}{6}\nonumber\\&+&\frac{ a(B+d_1)}{2 \beta_{10}} (T^l(\alpha))^2 +\frac{ a b}{4\beta^2_{10}} \frac{(T^l(\alpha))^4}{8}\bigg)  
\end{eqnarray}
and profit functional $J$ takes the form $J^r_{max} (\alpha)$
 \begin{eqnarray}
 &=&p(d_1 T^r(\alpha) +d_2 \frac{(T^r(\alpha))^2}{2}+d_3 \frac{(T^r(\alpha))^3}{3})-\bigg(c_d(l)+\frac{s_1}{T^r(\alpha)}\bigg)T^r(\alpha)\nonumber\\
 &-&a\bigg( C T^r(\alpha)+D \frac{(T^r(\alpha))^2}{2} -(2d_3 -\frac{ b }{2\beta_{10}})\frac{(T^r(\alpha))^4}{24}+(\frac{ a}{2\beta_{10}}-d_2) \frac{(T^r(\alpha))^3}{6}
\bigg)\nonumber\\
 &-&b\bigg( C \frac{T^r(\alpha)}{2}+D \frac{(T^r(\alpha))^3}{3} -(2d_3 -\frac{ b }{2\beta_{10}})\frac{(T^r(\alpha))^5}{30}+(\frac{ a}{2\beta_{10}}-d_2) \frac{(T^r(\alpha))^4}{8}\bigg)\nonumber\\
 &-&\bigg(c_{10}+\frac{s_2}{T^r(\alpha)}\bigg)\bigg( D T^r(\alpha)+ d_1 T^r(\alpha) +\frac{ b }{2\beta_{10}}\frac{(T^r(\alpha))^3}{6}+\frac{ a}{2\beta_{10}} \frac{(T^r(\alpha))^2}{2}\bigg) \nonumber\\
& -&\beta_{10}\bigg( (D+ d_1)^2 +(\frac{ b }{2\beta_{10}})^2\frac{(T^r(\alpha))^5}{20}
 +(\frac{ a}{2\beta_{10}})^2 \frac{(T^r(\alpha))^3}{3}+\frac{ b(D+d_1) }{\beta_{10}}\frac{(T^r(\alpha))^3}{6}\nonumber\\&+&+\frac{ a(D+d_1)}{2 \beta_{10}} (T^r(\alpha))^2 +\frac{ a b}{4\beta^2_{10}} \frac{(T^r(\alpha))^4}{8}\bigg)  
\end{eqnarray}
\section{Numerical Experiments}
\noindent To illustrate the models, we assume the following input data for the proposed Model. 
\subsection{Input: }The inputs are given in the following Table-1 
\begin{table}[h!]
\begin{center} \caption{Input values of Model-6} \centerline{}
\begin{tabular}{|cc c c c ccc c |} \hline
$L$   &$ N$&$  cd(l)$&$ h(t)$ & $ c_{10}$& $ \beta_{10} $&$p \mbox{(in \$)}$& $ s_1$ &
\\\hline
$ 40$ & $60$    & $100$     & $3+0.2 t^2$ 
&$0.7$ & $0.5$& $200$&  $ 10$ &\\\hline
$a$   &$ b $&$ b_1 $&$ d_1$ & $ d_2$& $ d_3 $&$T$& $ s_2$ &
\\\hline
$ 3$ & $0.2$    & $0.2$     & $7$ 
&$4$ & $2$& $12$&  $ 3$ &\\\hline
\end{tabular}
\end{center}
\end{table}
\subsection{Output:}
\noindent Using the above input data, the profit function of Model  has been  maximized and
with the help of GRG (LINGO-11.0), the unknowns, $x(t),u(t)$ and maximum profit etc. are evaluated. Here maximum profit is $\$247007.20$ for $\alpha=1$;  $\$305109.50$ for the right cut of $\alpha=0.4$ and $\$193834.20$ for the left cut of $\alpha=0.4$. The numerical values of $ u^r(t,\alpha), x^r(t,\alpha);u^l(t,\alpha),x^l(t,\alpha); d(t)$ for the Model are given in Table-2 and 3 respectively.  These are also graphically depicted in Figures-2 and 3 respectively. One comparison results for production and  stock level for right cut and left cut of $\alpha =0.4 $ are shown in Figure-4.   
\begin{table}[h!]
\begin{center}\scriptsize{} \caption{Values of $ u^r(t,\alpha), d(t), x^r(t,\alpha)$
for right cut of $\alpha=0.4$} \centerline{}
\begin{tabular}{|cc c c c  |} \hline
$\mbox{Time}$&$ u^r(t,\alpha)$&$  d({t})$&$x^r(t,\alpha)$ &
\\\hline
$ t=0$ & $123.95$    & $7.00$     & $0$ 
&  \\
$ t=1$ & $127.05$    & $13.00$     & $115.81$ 
&  \\
$ t=2$ & $130.35$    & $23.00$     & $226.83$ 
&  \\
$ t=3$ & $133.85$    & $37.00$     & $329.25$ 
&  \\
$ t=4$ & $137.55$    & $55.00$     & $419.27$ 
&  \\
$ t=5$ & $141.45$    & $77.00$     & $493.09$ 
&  \\
$ t=6$ & $145.55$    & $103.00$     & $546.91$ 
&  \\
$ t=7$ & $149.85$    & $133.00$     & $576.93$ 
&  \\
$ t=8$ & $154.35$    & $167.00$     & $579.34$ 
&  \\
$ t=9$ & $159.05$    & $205.00$     & $550.36$ 
&  \\
$ t=10$ & $163.95$    & $247.00$     & $486.18$ 
&  \\
$ t=11$ & $169.05$    & $293.00$     & $383.00$ 
&  \\
$ t=12$ & $174.35$    & $343.00$     & $237.02$ 
&  \\\hline
\end{tabular}
\end{center}
\end{table}
\begin{table}[h!]
\begin{center} \scriptsize{} \caption{Values of $ u^l(t,\alpha), d(t), x^l(t,\alpha)$
for left cut of $\alpha=0.4$ } \centerline{}
\begin{tabular}{|cc c c c  |} \hline
$\mbox{Time}$&$ u^l(t,\alpha)$&$  d({t})$&$x^l(t,\alpha)$ &
\\\hline
$ t=0$ & $86.27$    & $7.00$     & $0$ 
&  \\
$ t=1$ & $89.37$    & $13.00$     & $78.14$ 
&  \\
$ t=2$ & $92.67$    & $23.00$     & $151.47$ 
&  \\
$ t=3$ & $96.17$    & $37.00$     & $216.21$ 
&  \\
$ t=4$ & $99.87$    & $55.00$     & $268.55$ 
&  \\
$ t=5$ & $103.77$    & $77.00$     & $304.69$ 
&  \\
$ t=6$ & $107.87$    & $103.00$     & $320.83$ 
&  \\
$ t=7$ & $112.17$    & $133.00$     & $313.17$ 
&  \\
$ t=8$ & $116.67$    & $167.00$     & $277.90$ 
&  \\
$ t=9$ & $121.37$    & $205.00$     & $211.37$ 
&  \\
$ t=10$ & $126.27$    & $247.00$     & $109.38$ 
&  \\
$ t=11$ & $-$    & $-$     & $\mbox{Infeasible}$ 
&  \\
$ t=12$ & $-$    & $-$     & $\mbox{Infeasible}$ 
&  \\\hline
\end{tabular}
\end{center}
\end{table}
\begin{figure}[h!]
\centering
\includegraphics[width=6cm, height=5cm]{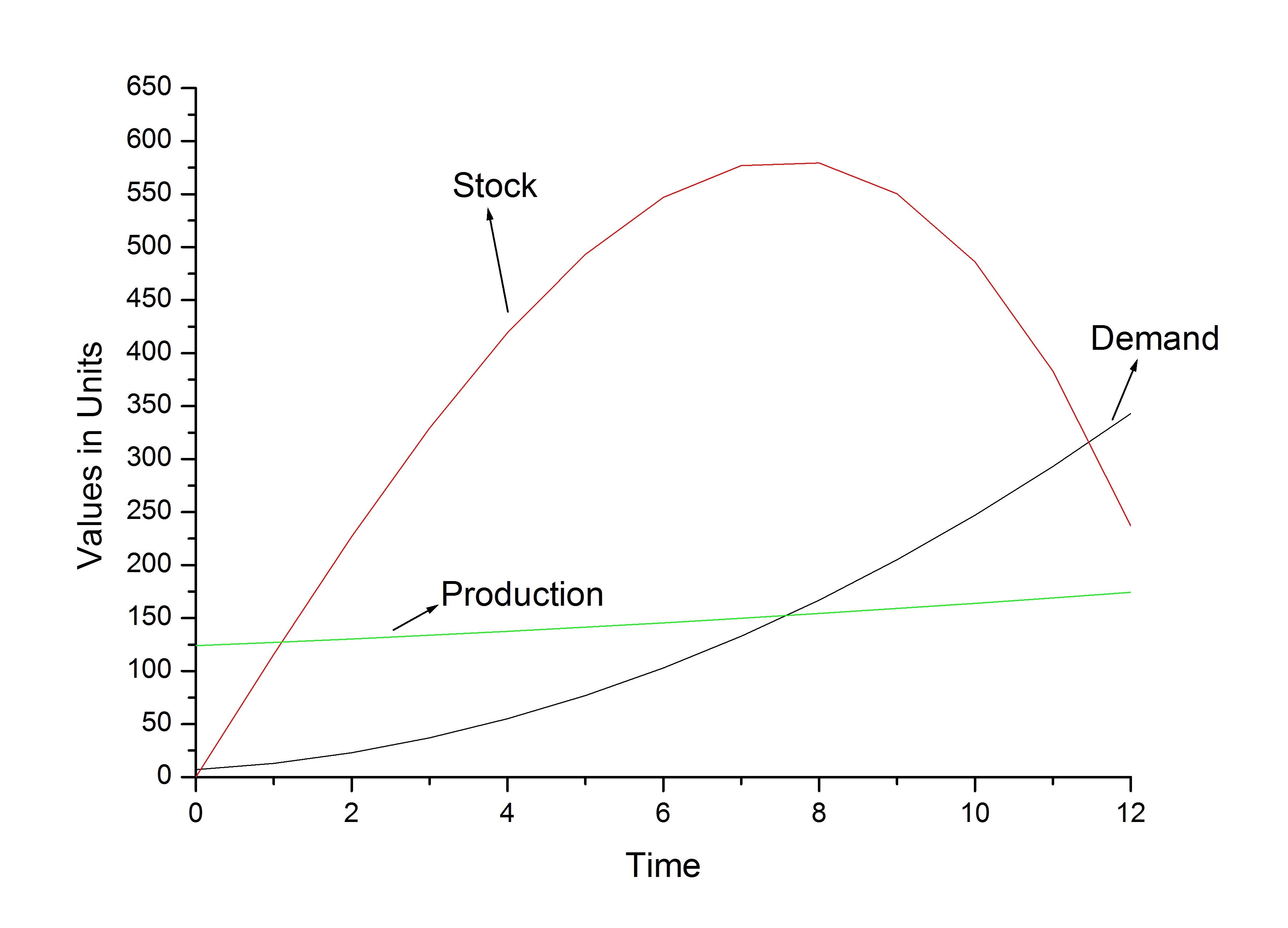}
\caption{Time Vs Production, Stock Level and Demand for Right cut of $\alpha$}
\end{figure}
\begin{figure}[h!]
\centering
\includegraphics[width=6.9 cm, height=5 cm]{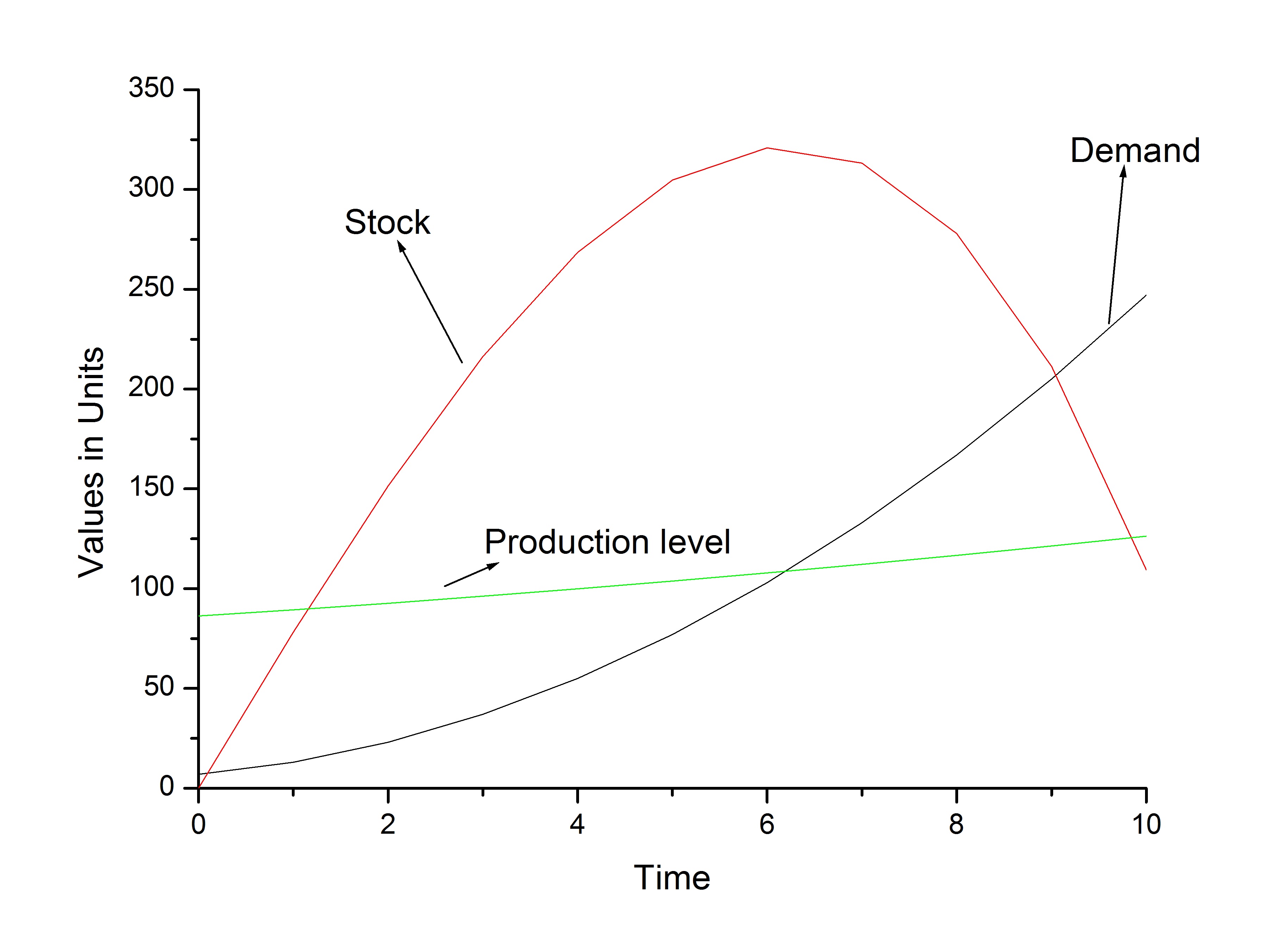}
\caption{Time Vs Production, Stock Level and Demand for Left cut of $\alpha$.}
\end{figure}
\begin{figure}[h!]
\centering
\includegraphics[width=7.5cm, height=6.75 cm]{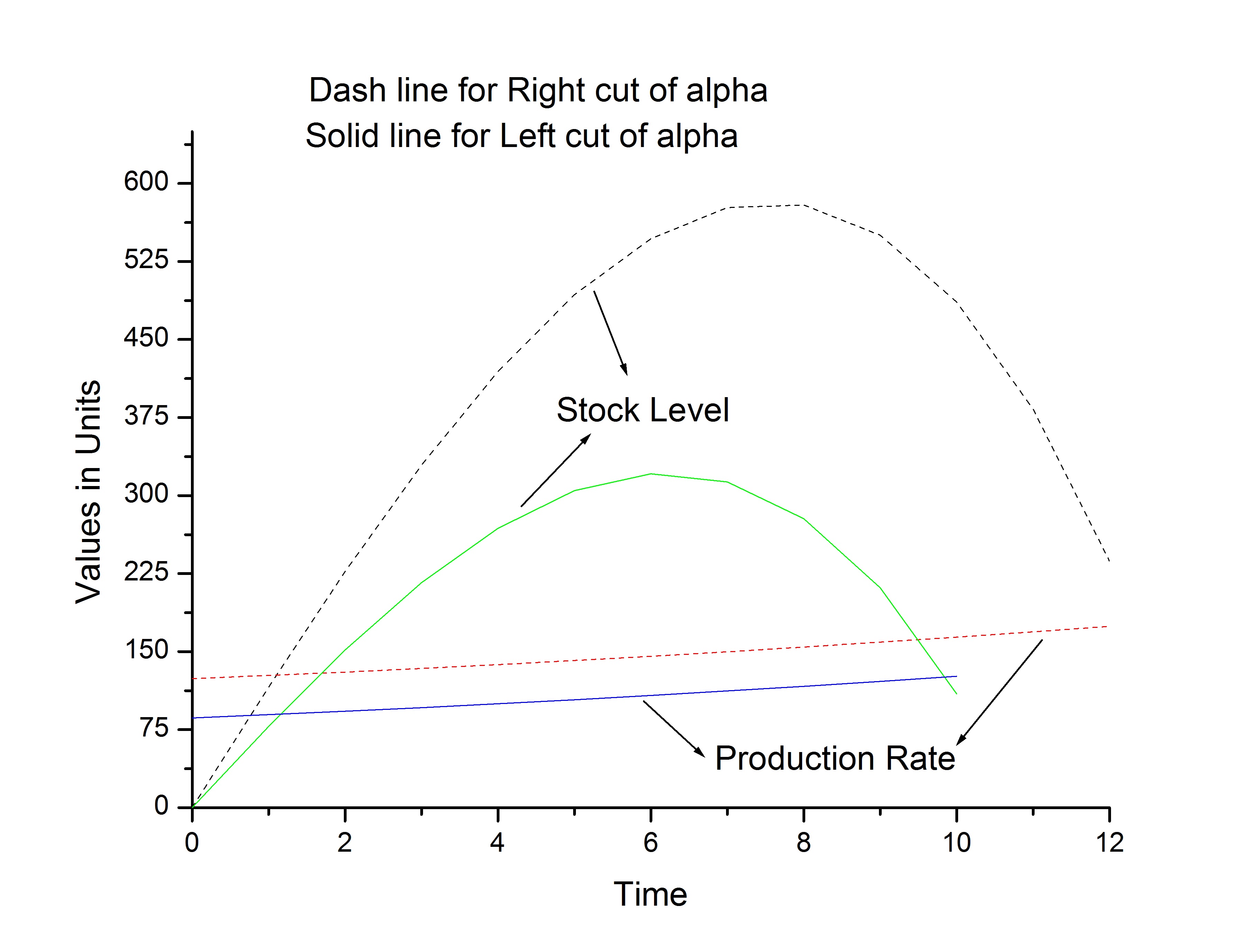}
\caption{Comparison of Time Vs Production  and Stock for Right cut and Left of $\alpha$ }
\end{figure}
\section{Sensitivity Analysis}
 For the  value of $\alpha=0.6$, the profit for left cut of $\alpha$ is $\$$ 210921.90  and for right cut of $\alpha $ is  $\$$ 285294.90   respectively. Changing the value of $\alpha$ as $\alpha=0.8$, the profit for left cut of $\alpha$ is  $\$$ 228661.60  and for right cut of $\alpha $ is $\$$ 265904.80  respectively. For those value of $\alpha $ (for right cut of $\alpha$) and with the help of GRG (LINGO-11.0) the unknowns,
$u^r(t,0.6), u^r(t,0.8); d(t)$ and $x^r(t,0.6), x^r(t,0.8)$;  for $\alpha=0.6,0.8 $ are evaluated which are given in Table-4.  Similarly, for left cut of $\alpha $, values are shown in Table-5.  Taking different values of $\alpha$, the profit values for right cut of $\alpha$ are given in Table-6 and for left cut of $\alpha$  are given in Table-7. Taking the value of $\alpha $ as $\alpha=0.6$  and $\alpha=0.8$ for right cut of $\alpha$ and  for left cut of $\alpha $, the comparison  diagram for production and stock level has been drawn in Figures-5 and 6 respectively. 
\begin{table}[h!]
\begin{center}\scriptsize 
\caption{Values of $u^r(t,0.6), u^r(t,0.8); d(t)$ and $x^r(t,0.6), x^r(t,0.8)$;  for $\alpha=0.6,0.8 $} 
\begin{tabular}{|c|c c c c c c c c c c c c c|} \hline
$t$   & $0$ & $1$ & $2$ & $3$ & $4$ &$5$ & $6$ & $7$& $8$& $9$ & $10$ & $11$ & $12$
\\\hline
$u^r(t,0.6)$ & $117.16$ & $120.26$ & $123.56$ &$127.06$ & $130.76$& $134.66$ &$138.76$ & $143.06$& $147.56$  & $152.26$& $157.16$ &$162.26$  & $167.56$ \\
$x^r(t,0.6)$ & $0.00$ & $109.11$ & $213.26$ &$308.89$ & $392.13$& $459.16$ &$506.19$ & $529.42$& $525.05$  & $489.28$& $418.32$ &$308.35$  & $155.58$
\\\hline
 $d(t)$ & $7.00$ & $13.00$ & $23.00$ & $37.00$ & $55.00$ & $77.00$  & $103.00$  & $133.00$  & $167.00$    & $205.00$  &$247.00$& $293.00$ & $343.00$ 
\\\hline
 $u^r(t,0.8)$ & $110.58$ & $113.68$&  $116.98$ &$120.48$ & $124.18$ & $128.08$ & $132.18$ & $136.48$ & $140.98$  & $145.68$& $150.58$ &$155.68 $  & $160.98$\\
 $x^r(t,0.8)$ & $0$ & $102.44$&  $200.09$ &$289.14$ & $365.79$ & $426.24$ & $466.68$ & $483.33$ & $472.38$  & $430.03$& $352.48$ &$235.92 $  & $76.57 $
\\\hline
\end{tabular}
\end{center}
\end{table}
\begin{table}[h!]
\begin{center}\scriptsize 
\caption{Values of $u^l(t,0.6), u^l(t,0.8); d(t)$ and $x^l(t,0.6), x^l(t,0.8)$;  for $\alpha=0.6,0.8 $} 
\begin{tabular}{|c|c c c c c c c c c c c c c|} \hline
$t$   & $0$ & $1$ & $2$ & $3$ & $4$ &$5$ & $6$ & $7$& $8$& $9$ & $10$ & $11$ & $12$
\\\hline
$u^l(t,0.6)$ & $92.04$ & $95.14$ & $98.45$ &$101.94$ & $105.64$& $109.54$ &$113.65$ & $117.94$& $122.44$  & $127.14$& $132.04$ &$137.14$  & $-$ \\
$x^l(t,0.6)$ & $0.00$ & $83.92$ & $163.02$ &$233.53$ & $291.64$& $333.56$ &$355.47$ & $353.58$& $324.09$  & $263.20$& $167.12$ &$32.03$  & $\mbox{Infeasible}$
\\\hline
 $d(t)$ & $7.00$ & $13.00$ & $23.00$ & $37.00$ & $55.00$ & $77.00$  & $103.00$  & $133.00$  & $167.00$    & $205.00$  &$247.00$& $293.00$ & $-$ 
\\\hline
 $u^l(t,0.8)$ & $98.02$ & $101.12$&  $104.42$ &$107.92$ & $111.62$ & $115.52$ & $119.62$ & $123.92$ & $128.42$  & $133.12$& $138.02$ &$143.12$  & $- $\\
 $x^l(t,0.8)$ & $0.00$ & $89.88$&  $174.97$ &$251.46$ & $315.55$ & $363.44$ & $391.32$ & $395.41$ & $371.90$  & $316.99$& $226.88$ &$97.00 $  & $\mbox{Infeasible}$
\\\hline
\end{tabular}
\end{center}
\end{table}
\begin{table}[h!]
\begin{center}\scriptsize 
\caption{The values of right cut of $\alpha$ Vs Profit  }
\begin{tabular}{|c|cc c c c c ccc|} \hline
$\alpha$  &$0.1$ & ${{0.2}}$ & $0.3$ & $0.4$ & $0.5$ & $0.6$ &$0.7$&$0.8$ &$0.9$
\\\hline
$\mbox{Profit} \,\,\, J\,\,\, in \,\,\, $\$$$ &$335460.30$& ${325273.20}$ & $315152.70$ & $305109.50$ &$295153.70$ & $285294.90$& $275542.40$&$265904.80$ &$256390.40$
\\\hline
 \end{tabular}
\end{center}
\end{table}
\begin{table}[h!]
\begin{center}\scriptsize 
\caption{The values of left cut of  $\alpha$ Vs Profit  }
\begin{tabular}{|c|cc c c c c ccc|} \hline
$\alpha$  &$0.1$ & $0.2$ & $0.3$ & $0.4$ & $0.5$ & $0.6$ &$0.7$&$0.8$&$0.9$ 
\\\hline
$\mbox{Profit} \,\,\, J\,\,\, in \,\,\, $\$$$ &$169509.90$& $177438.60$ & $185547.50$ & $193834.21$ &$202293.90$ & $210921.90$& $219713.00$&$228661.65$ &$237761.74$
\\\hline
 \end{tabular}
\end{center}
\end{table}
\begin{figure}[h!]
\centering
\includegraphics[width=7 cm, height=6.5 cm]{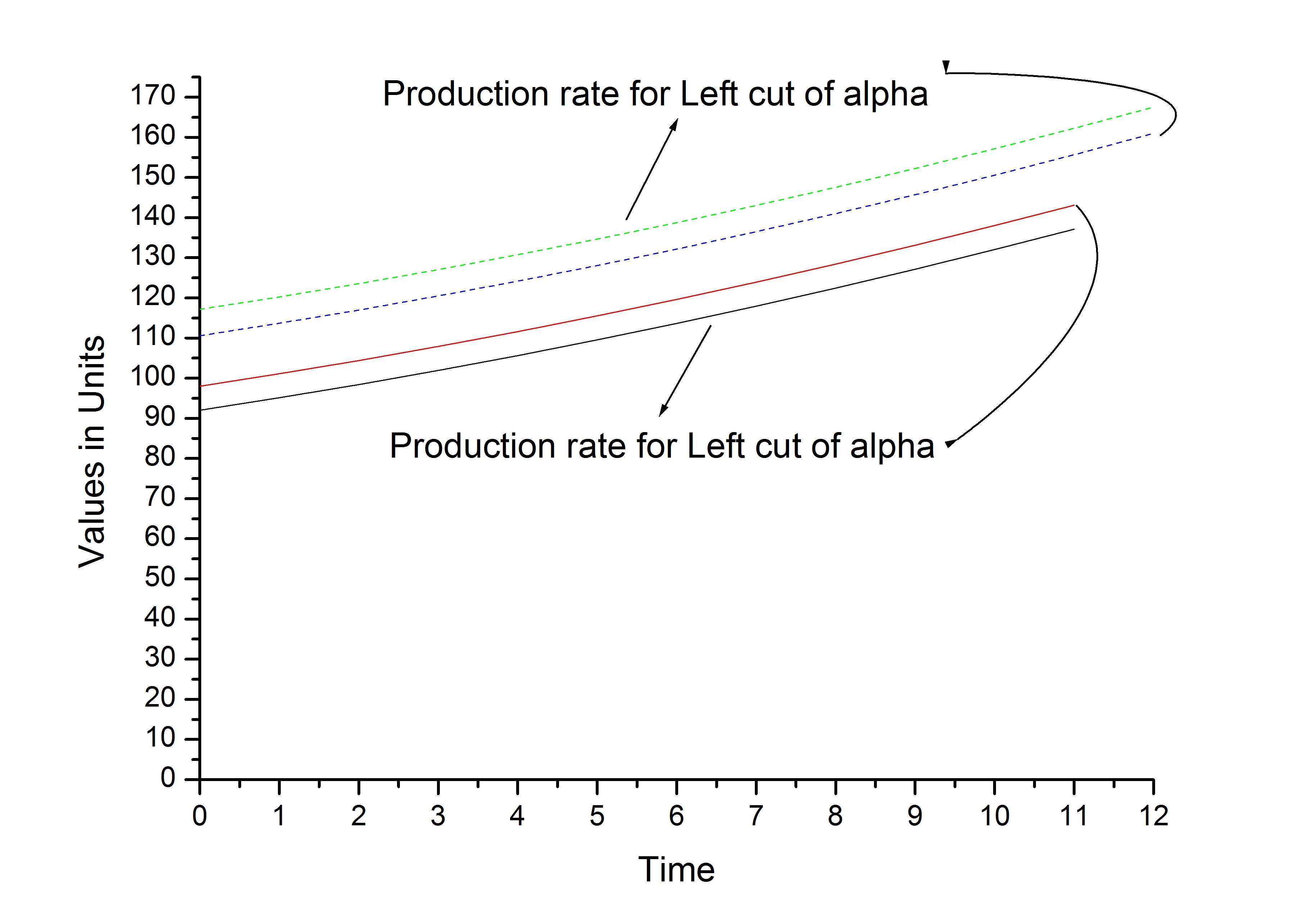}
\caption{Comparison of Time Vs Production for $ u^l(t,0.6), u^l(0.8)$  and $ u^r(t,0.6), u^r(0.8)$}
\end{figure}
\begin{figure}[h!]
\centering
\includegraphics[width=7 cm, height=6.5 cm]{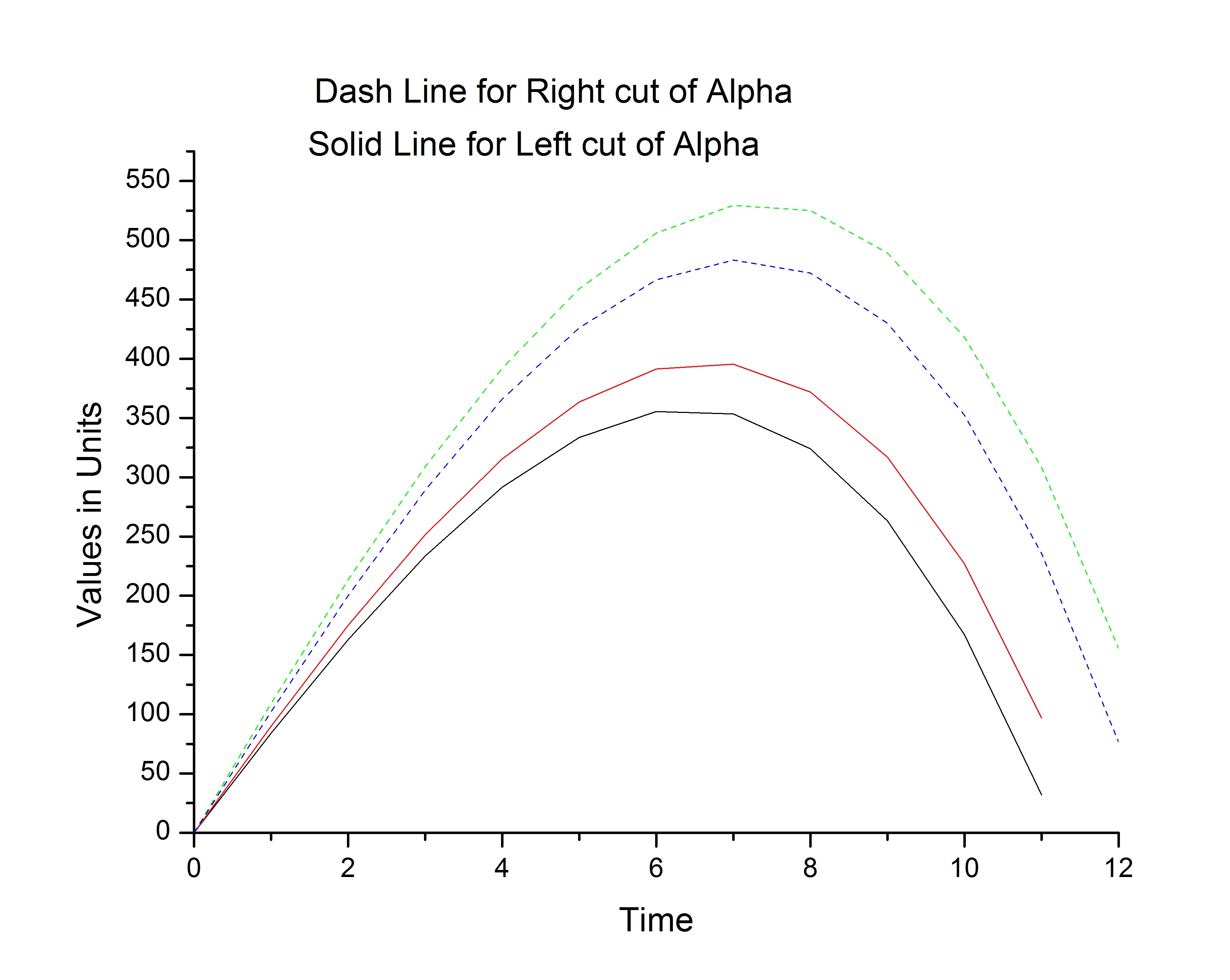}
\caption{Comparison of Time Vs Stock level for $ x^l(t,0.6), x^l(0.8)$  and $ x^r(t,0.6), x^r(0.8)$ }
\end{figure}
\newpage
\section{Discussion}
 In the Fig.s 2 and 3 for the models, initially production is more than demand and hence stock is built up after satisfying the demand. As demand is here dynamic which increases with time, after some time,  surpasses the product and then excess demand is satisfied from the stock. Again, the deficit in demand occurs, so the stock starts to decrease at this point and gradually decreases up to the end of the cycle
and finally becomes zero. 
\section{ Managerial insights}
In this investigation, a production-inventory model with fuzzy time period is formulated for maximum profit. Here, demand is of ramp type with respect to and production is a function of time. These considerations are appropriate for the production of seasonal products and fancy goods. Hence, the inventory  practitioners of these types of items may use this analysis to determine the   optimum production policy for maximum profit. Here, the end of season for seasonal products is not fixed, taken as imprecise. Demand of seasonal products such as winter garments, etc reaches maximum at a particular time of the period and then becomes zero at the end of the season. So the present investigation will be useful for the business people of seasonal products such as Festival greeting cards,winter garments, etc.

\section{Conclusion}
 Here a production-inventory model for the fuzzy time period is developed. For seasonality of an item, the demand should be adjusted before used in the forecast calculation of time. Seasonal products are sold only for a limited period of time every year. Actually, every year a seasonable product does not end at a specific time. In a year, several season products do have also different time periods. There is an intrinsic uncertainty in these time horizons. This uncertainty can be represented by fuzzy number. Thus the time periods of the items  made out of seasonable products are fuzzy. The models are also solved taking some of the inventory costs. The present models can be extended to the rough, fuzzy-rough, random, fuzzy-random environment taking constant part of holding cost, set-up cost, etc. as uncertain in nature.  The model can be extended to include multi-item fuzzy inventory problem with fuzzy space and budget constraints.


\begin{thebibliography}{0}
\bibitem{benjaffar2010} S. Benjaafar, J. P. Gayon and S. Tepex, Optimal control of a production-inventory system with customer impatience, Operations Research Letters 38 (4) (2010) 267-272
\bibitem{farhadinia2011} B. Farhadinia, Necessary optimality conditions for fuzzy variational problems, Information Sciences 181 (2011),1348-1357.
\bibitem{gabrielragsdell1977}G.A. Gabriel and K. M. Ragsdell, “The Generalized Reduced Gradient Method,” AMSE Journal of Engineering for Industry 99(1977), 384-400.
\bibitem {giri1996} B.C. Giri, A. Goswami and K.S. Chaudhuri, An EOQ model for deteriorating
items with time-varying demand and costs. Journal of the Operational Research Society 47 (11)(1996), 1398-1405.
\bibitem{hadleywhitin1975} G. Hadley and T.M. Whitin, An Optimal Final Inventory Model, Prentice Hall, 1975.
\bibitem{harris1913}F. Harris, Operations and Cost Factory Management Services, A.W. Shaw Co., Chicago, 1913.
\bibitem{hudong1995} J.Q. Hu and X. Dong, Optimal control for systems with deterministic production cycles, IEEE Trans. Autom. Contr. 40 (1995) 782-787
\bibitem{huloulon1995}J.Q. Hu and R. Loulon, Multi-product production/inventory control under random demands, IEEE Trans. Autom. Contr. 40 (1995) 350-355.
\bibitem{khouja1995} M. Khouja,  The Economic Production Lot-size Model Under Volume Flexibility. Computer and Oprations Research, 22(1995), 515-525.
\bibitem{kirk1970}D.E. Kirk, Optimal Control Theory. An Introduction, Prentice Hall, New Jersey, 1970.
\bibitem{maity2007}K. Maity, and M. Maiti,  “Possibility and Necessity Constraints and their Defuzzification – a Multi-item Production-inventory Scenario Via Optimal Control Theory,” European Journal of Operational Research 177(2007) 882-896.
\bibitem{naddor1966}E. Naddor, Inventory Systems, John Wiley, New York, 1966.
\bibitem{naidu2000}D.S. Naidu, Optimal Control System.(2000) Pocatello, ID: CRC Press.
\bibitem{najariyanfarahi2013} M. Najariyan and M.H. Farahi, Optimal control of fuzzy linear
controlled system with fuzzy initial conditions, Iranian
Journal of Fuzzy Systems 10(3) (2013), 21-35.
\bibitem{panda2008}D.Panda, S. Kar, K. Maity, and M. Maiti,  “A Single Period Inventory Model with Imperfect Production and Stochastic Demand Under Chance and Imprecise Constraints,” European Journal of Operational Research 188(2008) 121-139.
\bibitem{roul2015}  J.N.Roul, K. Maity, S. Kar,  and M. Maiti, Multi-item reliability dependent imperfect production inventory optimal control models with dynamic demand under uncertain resource constraint, International Journal of Production Research 53 (2015) 4993-5016.
\bibitem{roul2017}J. N. Roul, K. Maity, S. Kar,  and M. Maiti, Optimal control problem for an imperfect production process using fuzzy variational principle, Journal of Intelligent and Fuzzy Systems 32(1)(2017), 565-577.
\bibitem{salameh2000}M.K. Salameh,  and M. Y. Jaber, “Economic Production Quantity Model for Items with Imperfect Quality,” International Journal of Production Economics 64(2000) 59-64.
\bibitem{sarkar2011}B.Sarkar, S. S. Sana, and K. Chaudhuri, “An Imperfect Production Process for Time Varying Demand with Inflation and Time Value of Money – An EMQ Model,” Expert Systems with Applications 38(2011)13543-13548.
\bibitem{zequeria2004} R.I. Zequeira, B. Prida, and J. E. Valds, “Optimal Buffer Inventory and Preventive Maintenance for an Imperfect Production Process,” International Journal of Production Research 42 (5)(2004) 959-974.
\end{thebibliography}
\end{document}